\documentclass[english]{article}
\usepackage{amsmath}
\usepackage{amssymb}

\newcommand{\qed}{\hfill{\mbox{$\blacksquare$}}\medskip}
\newcommand{\proof}{\noindent{\sc Proof.\ }}

\newtheorem{definition}{Definition}[section]

\newtheorem{theorem}[definition]{Theorem}
\newtheorem{proposition}[definition]{Proposition}

\newtheorem{remark}[definition]{Remark}

\newtheorem{example}[definition]{Example}

\newcommand{\N}{\ensuremath{\mathbb{N}}}
\newcommand{\Z}{\ensuremath{\mathbb{Z}}}
\newcommand{\R}{\ensuremath{\mathbb{R}}}
\newcommand{\C}{\ensuremath{\mathbb{C}}}

\title{Behavior of bivariate interpolation operators at points of discontinuity of the first kind}
\author{Michele Campiti, Giusy Mazzone, Cristian Tacelli\thanks{Department of Mathematics ``E. De Giorgi, University of
Salento, P.O.Box 193, 73100 Lecce, Italy. E-mails: michele.campiti@unisalento.it, giusy.mazzone@unisalento.it, cristian.tacelli@unisalento.it.''}}
\date{}

\makeatother

\begin{document}

\maketitle

\pagestyle{plain}

\bigskip{}

\noindent \textit{AMS Classification (2010): 41A05, 40A30}

\bigskip{}




\begin{abstract}
  We introduce an index of convergence for double sequences of real numbers. This index is used to describe
the behaviour of some bivariate interpolation sequences at points of discontinuity of the first kind.
We consider in particular the case of bivariate Lagrange and Shepard operators.
\end{abstract}

\section{An index of convergence for double sequences}

In this paper we consider a general index of convergence for multiple sequences of real numbers. This index turns out to be useful in the description of non converging sequences and in some cases it can give complete information of the behavior of these sequences. 
This is the case for example of some phenomena in interpolation theory where we have at our disposal some results on the failure of the convergence at points of discontinuity of the function but a complete behavior has not yet obtained. In particular we shall concentrate ourselves on double sequences of real numbers and on the Lagrange and Shepard operators in the bivariate case, where we shall be able to furnish complete information on their behavior at points of discontinuity of the first kind of a function in terms of the index of convergence.

We start with the definition of the index of convergence for multi-indexed sequences of real numbers, which generalizes in a natural
way that of index of convergence for a sequence of real numbers given in \cite{CMT}.

In general, if $K\subset \N^m$,\ $m\ge1$, the lower density and, respectively, the upper density of $K$ can be defined by
\[
\delta_-(K):=\liminf_{n\to+\infty}\frac{\left|K\cap\{\,1,\dots,n\,\}^m\right|}{n^m}\,,\quad
\delta_+(K):=\limsup_{n\to+\infty}\frac{\left|K\cap\{\,1,\dots,n\,\}^m\right|}{n^m}\,.
\]
In the case where $\delta_-(K)=\delta_+(K)$ the density of $K$ is defined as follows
\[
\delta(K):=\delta_-(K)=\delta_+(K)\;.
\]

The equalities $\delta_-(K)=1-\delta_+(K^c)$ and $\delta_+(K)=1-\delta_-(K^c)$
remain true and can be shown as in \cite{CMT}.

It follows the definition of index of convergence.

\begin{definition} \label{df:index}
Let $(x_{n_1,\dots,n_m})_{n_1,\dots,n_m\geq1}$ be a multi-index sequence of real numbers.
If $L\in\R$, the \emph{index of convergence of the sequence} $(x_{n_1,\dots,n_m})_{n_1,\dots,n_m\geq1}$
\emph{to} $L$ is defined by
\[
i\left( x_{n_1,\dots,n_m};L \right):=1-\sup_{\varepsilon>0}\delta_+(\{\,(n_1,\dots, n_m)\in \N^m \mid x_{n_1,\dots,n_m}\notin ]L-\varepsilon,L+\varepsilon[\, \, \})\;.
\]
Moreover, we also set
\begin{align*}
i\left( x_{n_1,\dots,n_m};+\infty \right)&:=1-\sup_{M\in\R}\delta_+\left(\left\{\,(n_1,\dots,n_m)\in \N^m\mid x_{n_1,\dots,n_m}\notin ]-\infty,M]\, \right\}\right)\;,
\\
i\left( x_{n_1,\dots,n_m};-\infty \right)&:=1-\sup_{M\in\R}\delta_+\left(\left\{\,(n_1,\dots,n_m)\in \N^m\mid x_{n_1,\dots,n_m}\notin[M,+\infty[\, \right\}\right)\;.
\end{align*}
Finally, we can also define the \emph{index
of convergence of} $(x_{n_1,\dots,n_m})_{n_1,\dots,n_m\geq1}$ \emph{to} a subset $A$ of $\R$ as follows
\[
i(x_{n_1,\dots,n_m},A):=1-\sup_{\varepsilon>0}\delta_{+}(\{\,(n_1,\dots,n_m)\in \N^m\mid x_{n_1,\dots,n_m}\notin A+B_\varepsilon\, \, \})\;,
\]
where $B_\varepsilon:=]-\varepsilon,\varepsilon[^m$.
\end{definition}

In the case $m=1$ we obtain exactly the index of convergence considered in \cite{CMT}. Since all results can be easily extended from double sequences to multi-index sequences, for the sake of simplicity in the sequel we shall consider only the case $m=2$ of double sequences of real numbers.

\begin{remark} \label{rm:chardelta-} \rm
If necessary, we shall use the following explicit expressions of the index of convergence of a double sequence $(x_{n,m})_{n\geq1}$
\begin{align*}
i(x_{n,m};L)&=1-\sup_{\varepsilon>0}\delta_{+}\left(\{\,(n,m)\in \N^2\mid x_{n,m}\notin ]L-\varepsilon,L+\varepsilon[\, \} \right)
\\
&=1+\inf_{\varepsilon>0}\left(-\delta_{+}\left(\{\,(n,m)\in \N^2\mid \, x_{n,m}\notin ]L-\varepsilon,L+\varepsilon[\, \}
\right)\right)
\\
&=\inf_{\varepsilon>0}\left(1-\delta_{+}\left(\{\,(n,m)\in \N^2\mid x_{n,m}\notin ]L-\varepsilon,L+\varepsilon[\, \}
\right)\right)
\\
&=\inf_{\varepsilon>0}\delta_{-}\left(\{\,(n,m)\in \N^2\mid x_{n,m}\in ]L-\varepsilon,L+\varepsilon[\, \}\right)\;,
\end{align*}
and if $A\subset\R$
\[
i(x_{n,m},A)=\inf_{\varepsilon>0}\delta_{-}(\{\,(n,m)\in \N^2\mid x_{n,m}\in A+B_\varepsilon\, \})\;.
\]
\end{remark}

\begin{example} \label{ex:cos_equi} \rm
As a simple example, we can take $x_{n,m}:=\cos n\pi/2\,\cos m\pi/2$. It is easy to recognize that
\[
i(x_{n,m};0)=\frac 34,\ \ i(x_{n,m};1)=\frac 18,\ \ i(x_{n,m};-1)=\frac 18.
\]

\end{example}

In the next proposition we point out some relations between the index of convergence and the density of a
suitable converging subsequences.

\begin{proposition}\label{prop:char_index_multi}
Let $(x_{n,m})_{n,m \geq1}$ be a double sequence of real numbers and $\sigma\in]0,1]$. Then
$i(x_{n,m},L)\geq \sigma$ if and only if there exists a subsequence
$\left(x_{k(n,m)}\right)_{n,m\geq1}$ converging to $L$ such that
\[
\delta_{-} \left(\{\,k(n,m)\mid n,m\in \N\, \} \right)\ge \sigma\;.
\]
\end{proposition}
\proof $\Rightarrow$)\quad
For every $k\ge1$, we consider the set $M_{1/k}:=\{\,(n,m)\in\N^2\mid |x_{n,m}-L|<1/k\, \}$.
{}From Remark \ref{rm:chardelta-}, for every $k\in \N$ there exists $\tilde \nu_k$ such that
\[
\frac{\left| M_{1/k}\cap \{\,1,2,\dots,j\, \}^2\right|}{j^2}\geq \sigma-\frac{1}{k}
\]
whenever $j>\tilde \nu_k$. At this point we define recursively the sequence $(\nu_k)_{k\ge1}$
by setting $\nu_1=\tilde \nu_1$ and $\nu_k=\max \{\,\tilde \nu_k,\nu_{k-1}+1\, \}$. We have
\begin{equation}
  \label{eq:density-n2}
\frac{\left| M_{1/k}\cap \{\,1,2,\dots,j\, \}^2\right|}{j^2}\geq \sigma-\frac{1}{k}
\text{ for all } j>\nu_k\;.
\end{equation}

Consider the set of integers
\[
K=\bigcup_{k\ge1} \left( M_{1/k}\cap \{\,1,2,\dots,\nu_{k+1}\, \}^2\right)
\]
and the subsequence $\left\{\,x_{n,m}\mid (n,m)\in K\, \right\}$.

For every $\varepsilon >0$, let $\ell\in \N$ such that $1/\ell\leq \varepsilon$.
Then for every $(n,m)\in K$ satisfying $n,m>\nu_ \ell$
we have
$(n,m)\in  \bigcup_{k\geq  \ell} \left( M_{1/k}\cap \{\,1,2,\dots,\nu_{k+1}\, \}^2\right)$
and hence $|x_{n,m}-L|< \frac{1}{\ell}\leq \varepsilon$. This shows that the subsequence
$\left\{\,x_{n,m}\mid (n,m)\in K\, \right\}$ converges to $L$.

On the other hand, for every $j>\nu_\ell$, there
exists $\tilde\ell\geq \ell$ such that $\nu_{\tilde\ell}<j\leq \nu _{\tilde\ell +1}$ and
 thanks to \eqref{eq:density-n2}
we have
\begin{align*}
\frac{|K\cap\{\,1,2,\dots,j\, \}^2| }{j^2}&\geq \frac{\left| M_{1/\tilde\ell}\cap \{\,1,2,\dots,\nu_{\tilde\ell+1} \, \}^2\cap \{\,1,2,\dots,j\, \}^2\right|}{j^2}\\
&=\frac{\left| M_{1/\tilde\ell}\cap \{\,1,2,\dots, j\, \}^2\right|}{j^2}\geq \sigma-\frac{1}{\tilde\ell}\geq \sigma-\frac{1}{\ell} \geq \sigma-\varepsilon
\end{align*}
that is
\[
\liminf_{n \to \infty}\frac{|K\cap\{\,1,2,\dots,j\, \}^2| }{j^2}\geq \sigma\;.
\]

\noindent $\Leftarrow$)\quad
We suppose that there exists
a subsequence $(x_{k(n,m)})_{n\geq1}$ converging to $L$
such that such that
$\delta_{-} \left(\{\,k(n,m)\mid n,m\in \N\, \} \right)\ge \sigma$.
For every $\varepsilon > 0$ there exists $\nu_\varepsilon\in \N$ such
that $|x_{k(n,m)}-L|<\varepsilon$ whenever $n,m\geq \nu_\varepsilon$. Hence
\begin{align*}
\delta_-(\{\,(n,m)\in \N^2\mid  |x_{n,m}-L|<\varepsilon\, \})&\geq
\delta_-(\{\,(n,m)\in \N^2\mid |x_{k(n,m)}-L|<\varepsilon\, \})
\\
&=\delta_-(\{\,k(n,m)\mid  n,m\geq \nu_\varepsilon\, \})\\
&=\delta_-(\{\,k(n,m))\mid n,m\in \N\, \})\ge
\sigma
\end{align*}
and therefore, from Remark \ref{rm:chardelta-}, we obtain $i(x_{n,m},L)\geq \sigma$.
\qed

\begin{proposition}\label{pr:sum_gen_multi}
Let $(x_{n,m})_{n,m \geq1}$ be a double
sequence of real numbers and $(A_j)_{j\ge1}$ a sequence of subsets of $\R$ such that
$\overline{A_{k}}\cap \overline{A_{j}}=\emptyset$ for all $k\neq j$. Then
\[
0\leq \sum^{+\infty}_{k=1}i(x_{n,m},A_{k})\leq 1.
\]

In particular, if $(L_k)_{k\geq1}$ is a sequence of distinct elements of $[-\infty,\infty]$ such that, for every
$m\geq 1$
\begin{align*}
i(x_{n,m};L_k)=\alpha_k\;,
\end{align*}
for some $\alpha_{k}\ge0$, then
\[
0\leq \sum^{+\infty}_{k=1}\alpha_{k}\leq 1.
\]
\end{proposition}

\proof
Let $N\ge1$; since $\overline{A_{k}}\cap \overline{A_{j}}=\emptyset$ whenever $k\neq j$, we can choose
$\varepsilon$ such that
\[
(A_{k}+B_\varepsilon)\cap (A_{j}+B_\varepsilon)=\emptyset
\]
for all $k,j=1,\dots,N$, $k\neq j$.

Now consider the set
\[
M^{(k)}_{\varepsilon}:=\{\,(n,m)\in \N^2\mid x_{n,m}\in A_{k}+B_\varepsilon\, \}
\]
and observe that $M^{(k)}_{\varepsilon}\cap M^{(j)}_{\varepsilon}=\emptyset$
whenever $k,j=1,\dots,N$, $k\neq j$. Then we can conclude that
\begin{align*}
0\leq &\sum^{N}_{k=1}i(x_{n,m},A_{k})\leq \sum^{N}_{k=1}\delta_{-}(\{\,(n,m)\in \N^2\mid x_{n,m}\in A_{k}+B_\varepsilon
\, \})
\\
&=\sum^{N}_{k=1}\liminf_{n\rightarrow \infty}\frac{|M^{(k)}_{\varepsilon}\cap\{\,1,\dots,n\, \}^2|}{n^2}
\leq
\liminf_{n\rightarrow \infty}\left( \sum^{N}_{k=1}\frac{|M^{(k)}_{\varepsilon}\cap\{\,1,\dots,n\, \}^2|}{n^2} \right)
\\
&=\liminf_{n\rightarrow \infty}\frac{\left|\bigcup^{N}_{k=1}M^{(k)}_{\varepsilon}\cap\{\,1,\dots,n\, \}^2\right|}{n^2}
=\delta_{-}\left(\bigcup^{N}_{k=1}M^{(k)}_{\varepsilon}\right)\leq 1
\end{align*}
\qed

\begin{remark} \label{rm:sum=1} \rm
Observe that if in the preceding proposition we have $\sum_{k=1}^{+\infty}\alpha_{k}=1$, then every
subsequence $\left(x_{k(n,m)}\right)_{n,m\geq1}$ of $(x_{n,m})_{n,m\geq1}$
which converges to a limit $L$ different from each $L_k,$\ $k\ge1$,
necessarily satisfies $\delta_{-}(\{\,k(n,m)\mid n,m\in \N\, \})=0$ and therefore $i(x_{n,m};L)=0$.

Indeed, if
a subsequence $\left(x_{k(n,m)}\right)_{n,m\geq1}$ of $(x_{n,m})_{n,m\geq1}$ exists such that
$\delta_{-}(\{\,k(n,m)\mid n,m\in \N\, \})=\alpha>0$,
then by Proposition \ref{prop:char_index_multi} we get $i(x_{n,m},L)\geq \alpha$ and therefore
\[
i(x_{n,m},L)+\sum^{\infty}_{k=1}i(x_{n,m},L_{k})
\geq
\alpha+\sum^{\infty}_{k=1}\alpha_{k}>1
\]
which contradicts Proposition \ref{pr:sum_gen_multi}.
\end{remark}

\begin{proposition}\label{pr:1-bivariate}
Let $\left(x_{n,m}\right)_{n,m\geq1}$
be a double sequence of real numbers and
$\left(y_{n}\right)_{n\geq1}$ a sequence of real numbers. If
there exists a subsequence $(k(m))_{m\ge1}$ such that
$\lim _{m\to \infty }x_{n,k(m)}=y_n$ uniformly with respect to $n$ and if $\delta\left\{\,k(m)\mid m\in \N \, \right\}=\alpha$ then
\[
i(x_{n,m};A)\geq \alpha\, i(y_n;A)
\]
for every $A$ subset of $\R$.
\end{proposition}

\proof
Let us consider $K_1,K_2\subset \N$; since
\begin{align*}
&\delta_-(\{\, (i,j)\in \N^2\mid  (i,j)\in  K_1\times K_2  \, \})
=\delta_-(\{\, (i,j)\in \N^2\mid  i\in K_1, j\in K_2  \, \})
\\
&\quad=
\liminf_{n\to \infty}
\frac{|\{\,(i,j)\in \N^2| i\in K_1, j\in K_2\, \}\cap \{\,1,\dots,n\, \}^2|}{n^2}\\
&\quad
=\liminf_{n\to \infty}
\frac{|\{\,i\in\N\mid i\in K_1\, \}\cap \{\,1,\dots,n\, \}|}{n}
\frac{|\{\,j\in\N\mid j\in K_2\, \}\cap \{\,1,\dots,n\, \}|}{n}\\
&\quad
\geq \liminf_{n\to \infty}
\frac{|\{\,i\in\N\mid i\in K_1\, \}\cap \{\,1,\dots,n\, \}|}{n}
\liminf_{n\to \infty}
\frac{|\{\,j\in\N\mid j\in K_2\, \}\cap \{\,1,\dots,n\, \}|}{n}\\
&\quad=\delta_-(\{\,i\in \N\mid i\in K_1 \, \})
\delta_-(\{\,j\in \N\mid j\in K_2 \, \})\;,
\end{align*}
we have
\begin{equation}\label{eq:delta-prod}
 \delta_-(K_1\times K_2)\geq\delta_-(K_1)\delta_-(K_2)\;.
\end{equation}
For every $\varepsilon>0$ there exists $\eta\in\N$ such that $|x_{n,k(m)}-y_n|<\varepsilon$ whenever $m\geq \eta$ and $n\in \N$.
Then
\[
\{\,(n,k(m))\in \N^2\mid m\geq \eta\;,\;y_n\in A+B_\varepsilon \, \}
\subset \{\,(n,m)\in\N^2\mid  x_{n,m}\in A+B_{2\varepsilon}\, \}
\]
and consequently
\begin{align*}
 &\delta_-\left( \{\,(n,m)\in\N^2\mid  x_{n,m}\in A+B_{2\varepsilon}\, \}\right)\geq
\delta_-\left( \{\,(n,k(m))\in \N^2\mid m\geq \eta\;,\;y_n\in A+B_\varepsilon \, \} \right)\\
&\quad \geq
\delta_-(\{\,n\in \N\mid y_n\in A+B_\varepsilon\, \})
\delta_-(\{\,k(m)\mid  m\geq \eta\, \})\\
&\quad =\alpha \delta_-(\{\,n\in \N\mid y_n\in A+B_\varepsilon\, \})
\geq \alpha i(y_n;A)\;.
\end{align*}
Taking the infimum with respect to $\varepsilon$ we obtain the desired result.
\qed

\begin{proposition}\label{pr:prod-succ}
Let $(x_n)_{n\geq 1}$ and $(y_n)_{n\geq 1}$ sequences of real numbers,
let $f:[0,1]\to [0,1]$ and $g:[0,1]\to [0,1]$ be injective
differentiable functions with differentiable inverses.
If $i(x_n;A)=|f^{-1}(A)|$ and $i(y_n;A)=|g^{-1}(A)|$ for every Peano-Jordan measurable set $A\subset [0,1]$,
then
\[
 i(x_ny_m,A)=|G^{-1}(A)|\;,
\]
where $G(x,y)=g(x)f(y)$.

\end{proposition}

\proof
Firstly we prove that $i(x_n;I)=\delta(\{\,n\in\N\mid x_n\in I\, \})$ for every interval $I\subset \R$.
Notice that for every interval $[a,b]\subset \R$
and for all $\varepsilon> 0$
\[
\delta_-(\{\,n\in \N|\, x_n \in ]a-\varepsilon,b+\varepsilon[\, \})\geq
\delta_-(\{\,n\in \N|\, x_n \in [a,b]\, \})
\]
and then, taking the infimum over $\varepsilon>0$,
\begin{equation}\label{eq:1g}
i(x_n;[a,b])\geq \delta_-(\{\,n\in \N|\, x_n\in [a,b]\, \}).
\end{equation}
On the other hand, for every $\delta >0$ we have that
\begin{align*}
&|f^{-1}([a+\delta,b-\delta])|=\inf_{\varepsilon>0}\delta_-(\{\,n\in\N\mid x_n\in]a+\delta-\varepsilon,b-\delta+\varepsilon[\, \})\\
&\quad \leq
\delta_-(\{\,n\in\N\mid x_n\in [a,b]\, \})\;.
\end{align*}
Since $f^{-1}$ is continuous, the function $\delta \mapsto |f^{-1}([a+\delta,b-\delta])|$ is continuous at $0$ and taking the limit as $\delta\rightarrow 0$ we have
\begin{equation*}
i(x_n;[a,b])=|f^{-1}([a,b])|\leq \delta_-(\{\,n\in \N|\, x_n\in [a,b]\, \})\;,
\end{equation*}
which jointly with \eqref{eq:1g} yields
\[
i(x_n;I)=\delta_-(\{\,n\in\N\mid x_n\in I\, \})
\]
for every interval $I$. Finally  we have
\begin{align*}
&\delta_+(\{\,n\in \N\mid  x_n\in I  \, \})=1-\delta_-(\{\,n\in \N\mid  x_n\in I^c  \, \})=1-|f^{-1}(I^c)|\\
&\quad =|f^{-1}(I)|=\delta_-(\{\,n\in \N\mid  x_n\in I  \, \})\;,
\end{align*}
so
\[
i(x_n;I)=|f^{-1}(I)|=\delta(\{\,n\in\N\mid x_n\in I\, \})\;.
\]

Let $I,J$ be real intervals, we have
\[
\delta(\{\, (n,m)\in \N\mid  (x_n,y_m)\in I\times J    \, \})
=\delta(\{\,n\in \N\mid x_n\in I \, \})
\delta(\{\,m\in \N\mid y_m\in J \, \})\;,
\]
indeed
\begin{align*}
&\delta(\{\, (i,j)\in \N^2\mid  (x_i,y_j)\in  I\times J  \, \})
=\delta(\{\, (i,j)\in \N^2\mid  x_i\in I, y_j\in J  \, \})
\\
&\quad=
\lim_{n\to \infty}
\frac{|\{\,(i,j)\in \N^2| x_i\in I, y_j\in J\, \}\cap \{\,1,\dots,n\, \}^2|}{n^2}\\
&\quad
=\lim_{n\to \infty}
\frac{|\{\,i\in\N\mid x_i\in I\, \}\cap \{\,1,\dots,n\, \}|}{n}
\frac{|\{\,j\in\N\mid y_j\in J\, \}\cap \{\,1,\dots,n\, \}|}{n}\\
&\quad
= \lim_{n\to \infty}
\frac{|\{\,i\in\N\mid x_i\in I\, \}\cap \{\,1,\dots,n\, \}|}{n}
\lim_{n\to \infty}
\frac{|\{\,j\in\N\mid y_j\in J\, \}\cap \{\,1,\dots,n\, \}|}{n}\\
&\quad=\delta(\{\,i\in \N\mid x_i\in I \, \})
\delta(\{\,j\in \N\mid y_j\in J \, \})\;.
\end{align*}

Then we have
\[
\delta\left( \{\, (n,m)\in\N^2\mid  (x_n,y_m)\in I\times J \, \}\right)
=|f^{-1}(I)||g^{-1}(J)|=|(f,g)^{-1}(I\times J)|;
\]
moreover, thank to the linearity of the limit, if $Q$ is a pluri-interval of $\R^2$ we have
\[
\delta\left( \{\, (n,m)\in\N^2\mid  (x_n,y_m)\in Q \, \}\right)
=|(f,g)^{-1}(Q)|\;.
\]

Now let $\tilde A$ be a Peano-Jordan measurable subset of $[0,1]^2$ and fix $\varepsilon>0$.
Since $F:=(f,g)^{-1}$ is a diffeomorphism from $[0,1]^2$ into $[0,1]^2$, the subset
$F(\tilde A)$ is measurable and therefore there exist pluri-intervals $Y_1,Y_2$ such
that $Y_1\subset F(\tilde A)\subset Y_1$ and $|Y_1|-|Y_2|\leq \varepsilon$.
Let  $Q_1=F^{-1}(Y_1)$ and $Y_2=F^{-1}(Q_2)$; these set are pluri-intervals and $Q_1\subset \tilde A\subset Q_2$.
We have
\begin{align*}
&\delta \left(  \{\, (n,m)\in\N^2\mid  (x_n,y_m)\in Q_1 \, \} \right)\leq
\delta \left(  \{\, (n,m)\in\N^2\mid  (x_n,y_m)\in \tilde A \, \} \right)\\
&\quad \leq
\delta \left(  \{\, (n,m)\in\N^2\mid  (x_n,y_m)\in Q_2 \, \} \right),
\end{align*}
from which
\begin{align*}
|Y_1|=|F(Q_1)|
\leq
\delta \left(  \{\, (n,m)\in\N^2\mid  (x_n,y_m)\in \tilde A \, \} \right)
\leq |F(Q_2)|=|Y_2|\;.
\end{align*}
On the other hand
\[
|Y_1|\leq |F(\tilde A)| \leq  |Y_2|
\]
and therefore
\begin{align*}
&\left|\delta \left(  \{\, (n,m)\in\N^2\mid  (x_n,y_m)\in \tilde A \, \} \right)
-|F(\tilde A)|\right|\leq |Y_2|-|Y_1|\leq \varepsilon\;
\end{align*}
and we get
\[
\delta
\left(
\{\, (n,m)\in\N^2\mid  (x_n,y_m)\in \tilde A \, \}
\right)
=|F(\tilde A)|\;.
\]

Finally we consider a Peano-Jordan measurable set $A\subset [0,1]$ and
the function $P:\R^2\to \R$, $P(x,y)=xy$. We have
\begin{align*}
\{\,(n,m)\in \N^2\mid  x_ny_m\in A \, \}&=
\{\,(n,m)\in \N^2\mid  P(x_n,y_m)\in A\, \}\\
&=\{\,(n,m)\in \N^2\mid  (x_n,y_m)\in P^{-1}(A)\, \}
\end{align*}
and hence
\begin{align*}
\delta( \{\, (n,m)\in \N^2\mid x_ny_m\in A\, \})&=
\delta\left( \{\, (n,m)\in \N^2\mid  (x_n,y_m)\in P^{-1}(A)\, \}\right)\\
&= |(f,g)^{-1}P^{-1}(A)| =|(f\cdot g)^{-1}(A)|
\end{align*}
where $(f\cdot g)(x,y)=G(x,y)=f(x)g(y)$.
Therefore for every Peano-Jordan measurable set $A\subset [0,1]$ we have
\begin{align}\label{eq:index_prod}
i(x_ny_m;A)&=\inf_{\varepsilon>0}\delta(\{\,(n,m)\in \N^2|\, x_ny_m\in A+B_\varepsilon\, \})
=\inf_{\varepsilon>0}|G^{-1}(A+B_\varepsilon)|
\nonumber
\\
&\geq |G^{-1}(A)|.
\end{align}
To prove the converse inequality, we argue by contradiction and suppose that
$i(x_ny_m;A)>|G^{-1}(A)|$. So there exists $\delta>0$ such that $i(x_ny_m;A)=|G^{-1}(A+B_\delta)|$.
Notice that $\overline{A}\cap \overline{(A+B_{\delta/2})^c}=\emptyset$ and, by Proposition
\ref{pr:sum_gen_multi}, we get
\[
i(x_ny_m;A)+i(x_ny_m;(A+B_{\delta/2})^c)\leq 1.
\]
Since $G^{-1}([0,1]^2)\subset [0,1]^2$, by \eqref{eq:index_prod} we have
\[
|G^{-1}(A+B_\delta)|\leq 1-|G^{-1}((A+B_{\delta/2})^c)|=|G^{-1}(A+B_\delta)|\;.
\]
This leads to a contradiction since the map
$\delta>0\rightarrow |G^{-1}(A+B_\delta)|$ is monotone increasing. Then our claim is achieved.
\qed

\begin{example}
As a further example, let $\alpha,\gamma \in [0,1)$ be irrational, $\beta,\delta \in [0,1)$ and consider
\[
x_{n,m}:=(n\alpha +\beta-[n\alpha+\beta])(m\gamma +\delta-[m\gamma +\delta])
\]
where $[x]$ denotes the integer part of $x$.

We already know that (see \cite[Example 1.4)-ii)]{CMT})
\[
i\left(n\alpha +\beta-[n\alpha+\beta];A\right)=|A|
\]
for every Peano-Jordan measurable set $A\subset [0,1[$, where $|\cdot|$ denotes the Peano-Jordan measure.

Now, applying Proposition \ref{pr:prod-succ} with the identity function in place of $f$ and $g$, we get
%
%
\[
i(x_{n,m};A)=\delta \left(\left\{\,(n,m)\in \N^2\mid x_{n,m}\in A\, \right\}\right)=|G^{-1}(A)|
\]
for every Peano-Jordan measurable set $A\subset [0,1[$, where $G(x,y)=x\,y$.
\end{example}

\section{Bivariate Lagrange operators on discontinuous functions}\label{sc:Lagrange}

We begin by considering the univariate Lagrange operators $(L_n)_{n\geq 1}$ at the Chebyshev nodes of second
type, which are defined by
\[
L_nf(x)=\sum_{k=1}^n\ell_{n,k}(x)f(x_{n,k}),
\]
where $f$ is a suitable function from $[-1,1]$ to $\R$,
\[
x_{n,k}=\cos\theta_{n,k},\qquad\theta_{n,k}=\frac{k-1}{n-1}\pi\qquad k=1,\ldots,n\;,
\]
are the Chebyshev nodes of second type and
\[
\ell_{n,k}(x)=\prod_{i\not=k}\frac{x-x_{n,i}}{x_{n,k}-x_{n,i}}
\]
are the corresponding fundamental polynomials.

Setting $x=\cos\theta$, with $\theta\in[0,\pi]$, the polynomials $\ell_{n,k}$ can be rewritten as follows
\[
\ell_{n,k}(\cos\theta)=\frac{(-1)^{k}}{(n-1)(1+\delta_{k,1}+\delta_{k,n})}\,\frac{\sin((n-1)\theta)\sin \theta}
{\cos\theta-\cos\theta_{n,k}}\;,
\]
where $\delta_{i,j}$ denotes the Kronecker symbol, that is
\begin{equation*}
\delta_{i,j}:=
\begin{cases}
0&\ \text{if }i\neq j,
\\
1&\ \text{if }i=j.
\end{cases}
\end{equation*}

Our first aim is to study the behavior of the sequence of Lagrange operators for a particular class of functions having a finite number of points of discontinuity of the first kind. This will simplify the subsequent discussion on the bivariate case.

We consider the function
$h_{x_0,d}:[-1,1]\rightarrow\R$ defined by
\begin{equation} \label{eq:def_h}
h_{x_0,d}(x):=\left\{\,
\begin{array}{ll}
0\;,\quad& x<x_0\;,\\
d\;,\quad& x=x_0\;,\\
1\;,\quad& x>x_0\;,
\end{array} \right. \qquad x\in[-1,1]\;,
\end{equation}
where $d$ is a fixed real number.

We also need to define the function $g:]0,1[\mapsto \R$ by setting
\begin{equation}\label{eq:g}
g(x):=
\frac{ \sin \left(\pi x\right) }{\pi}\,
J(1,x)\;,\qquad \text{if } x\in ]0,1[\;,
\end{equation}
where $J(s,a)$ denotes the Lerch zeta function
\[
J(s,a):=\sum^{+\infty}_{n=0}\frac{(-1)^n}{(n+a)^s}\;,\qquad a\in ]0,1]\;,\quad \Re [s]>0\;.
\]

The following result describes the behavior of Lagrange operators at the point $x_0$ in terms of the index of convergence defined in \cite{CMT} and corresponding to the case $m=1$ in Definition \ref{df:index}.

\begin{theorem} \label{th:Lagrange}
Let $x_0=\cos\theta_0\in]-1,1[$ and consider the functions $h:=h_{x_0,d}$ defined by \eqref{eq:def_h}. Then, the sequence of functions $\left( L_nh\right)_{n\ge1}$ converges uniformly to $h$
on every compact subsets of $[-1,1]\setminus\{\,x_0\, \}$.

As regards the behaviour of the sequence $(L_nh(x_0))_{n\ge1}$ we have the following cases.
\begin{itemize}
\item[i)] If $\frac{\theta_0}{\pi}=\frac pq $ with $p,q\in \N$, $q\neq 0$ and $GCD(p,q)=1$, then
  \[
    i\left(L_nh(x_0);d\right)=\frac1q\;,\qquad
    i\left(L_nh(x_0);
    g\left( \frac{m}{q}\right)\right)=\frac1q\;,
    \quad m=1,\dots, q-1.
  \]

\item [ii)] If $\frac{\theta_0}{\pi}$ is irrational and if
$A\subset \R$ is a Peano-Jordan measurable set,
then
\[
i\left(L_nh(x_0);A\right)=|g^{-1}(A)|\;,
\]
where $|\cdot|$ denotes the Peano-Jordan measure.
\end{itemize}
\end{theorem}

\proof
Let $a=\cos \theta_1\in [-1,x_0[$ and
$x=\cos \theta\in [-1,a]$;
for sufficiently large $n\ge1$ there exists  $k_0$ such that $0\leq\theta_{n,k_0}\leq \theta_0<\theta_{n,k_0+1}
<\theta_1\leq \theta\leq \pi$ and therefore
\[
0<\cos \theta_0-\cos \theta_1\leq \cos \theta_{n,k_0}-\cos \theta\;.
\]

We have $L_nh(\cos \theta)
=\sum_{k=1}^{k_0-1}
  \ell_{n,k}(\cos \theta)
  +d\ell_{n,k_0}(\cos \theta)$ if $\theta_{n,k_0}=\theta_0$, and
$L_nh(\cos \theta)=\sum_{k=1}^{k_0}
  \ell_{n,k}(\cos \theta)
$ if $\theta_{n,k_0}<\theta_0$;
hence
\begin{align*}
&L_nh(\cos \theta)\\&\quad=
  \sum^{k_0}_{k=1}\frac{(-1)^{k}}{(n-1)(1+\delta_{k,1})}\,\frac{\sin((n-1)\theta)\sin \theta}
  {\cos\theta-\cos\theta_{n,k}}
+(d-1)\chi_{\left\{\theta_{n,k_0}\right\}}(\theta_0)\ell_{n,k_0}(\cos \theta)
  \\
&\quad=\frac{\sin ((n-1)\theta)\sin\theta}{2(n-1)(\cos\theta-1)}
  +\sum^{k_0}_{k=1}\frac{(-1)^{k-1}}{n-1}\,\frac{\sin((n-1)\theta)\sin \theta}
  {\cos\theta_{n,k}-\cos\theta}\\
&\qquad+(d-1)\chi_{\left\{\theta_{n,k_0}\right\}}(\theta_0)\ell_{n,k_0}(\cos \theta)\;.
\end{align*}

The function $t\to \frac{1}{\cos t-\cos \theta}$
is positive and monotone increasing on the interval $[0, \theta[$;
since $0<\theta_{n,k}<\theta_{n,k+1}<\theta$ for every  $1\leq k\leq k_0$, we have
\begin{align*}
|L_n h(x)|&=|L_n(h)(\cos \theta)|\\
&\leq \frac{1}{2(n-1)(1-\cos\theta)}
  +\left| \frac{\sin((n-1)\theta)\sin \theta}{n-1}\,
   \frac{1}{\cos\theta_{n,k_0}-\cos\theta}\right| \\
&\qquad+|d-1| \left| \frac{\sin((n-1)\theta)\sin \theta}{n-1}\,
   \frac{1}{\cos\theta_{n,k_0}-\cos\theta}\right|
\\
  &\leq
\frac{1}{2(n-1)(1-\cos\theta_1)}+
\frac{1+|d-1|}{n-1}\frac{1}{\cos\theta_{0}-\cos\theta_1}.
\end{align*}
It follows that $\left(L_n h\right)_{n\ge1}$ converges uniformly to $h$ in $[-1,a]$.

Now let
$b=\cos \theta_2\in ]x_0,1[$ and
$x=\cos \theta\in [b,1]$. For sufficiently large $n\ge1$   there
exists $k_0$ such that $0\leq \theta\leq \theta_2<\theta_{n,k_0}
\leq \theta_0<\theta_{n,k_0+1}\leq 2\pi$
and consequently
\[
0<\cos \theta_2-\cos\theta_0\leq  \cos \theta-\cos\theta_{n,k_0+1}\;.
\]
Then
\begin{align*}
|1-&L_nh(x)|=|1-L_nh(\cos \theta)|=\left| \sum_{k=1}^{n}\ell_{n,k}(\cos \theta)-
\sum_{k=1}^{k_0}\ell_{n,k}(\cos \theta)h(\cos \theta_{n,k})\right|\\
&=\left|
  \frac{(-1)^{n+1}}{2(n-1)}\,
  \frac{\sin ((n-1)\theta)\sin \theta}{\cos\theta+1}\right.\\
&\quad\left.
  +\sum_{k=k_0+1}^{n}
  \frac{(-1)^{k}}{n-1}\,
  \frac{\sin ((n-1)\theta)\sin \theta}{\cos\theta-\cos\theta_{n,k}}\,-
  (d-1)\chi_{\left\{\theta_{n,k_0}\right\}}(\theta_0)\ell_{n,k_0}(\cos \theta)\right|\\
&\leq \frac{1}{n-1}\left[ \frac{1}{2(\cos\theta +1)}+
  \left|
  \frac{1} {\cos \theta- \cos \theta _{n,k_0+1}}
  \right|
    +|d-1|\left|
    \frac{1} {\cos \theta- \cos \theta _{n,k_0}}
  \right|\right]
\\
&\leq \frac{1}{2(n-1)(\cos\theta_2+1)}+
  \frac{1+|d-1|}{n-1}\frac{1}{\cos\theta_2-\cos\theta_{0}}\;,
\end{align*}
since the function
$t\to \frac{1}{\cos \theta-\cos t}$ is positive and monotone decreasing in
$]\theta,\pi]$ and
$\theta<\theta_{n,k-1}<\theta_{n,k}<\pi$ for every $k_0+1\leq k\leq n$.
So $\left(L_n h\right)_{n\ge1}$ converges uniformly to $h$ also in $[b,1]$.

Now, we study the behavior of $(L_nh(x_0))_{n\ge1}$.

We identify $x_0=\cos\theta_0$, $\theta_0\in [0,\pi]$.
For sufficiently large $n\ge1$ there exists  $k_0$ such that
$\theta_{n,k_0}\leq \theta_0<\theta_{n,k_0+1}$.
Let us denote $\displaystyle   \sigma_n=\frac{n-1}{\pi}(\theta_0-\theta_{n,k_0})$.
{}From $\frac{k_0-1}{n-1}\pi \leq\theta_0<\frac{k_0}{n-1}\pi$ we have that
$0\leq \sigma_n <1$;
then
\[
n-1=\frac {\pi}{\theta_0}(\sigma_n+k_0-1)
\]
and moreover
\[
k_0-1\leq \frac{n-1}{\pi}\theta_0< k_0\;,
\]
that is $k_0-1=\left[\frac{n-1}{\pi}\theta_0\right]$
and
\begin{equation} \label{eq:sigman}
\sigma_n=\frac{n-1}{\pi}\theta_0-\left[ \frac{n-1}{\pi}\theta_0\right]\;.
\end{equation}

If $x_0$ is a Chebyshev node, that is $\theta_0=\theta_{n,k_0}$ and $\sigma_n=0$, then
\begin{equation}\label{eq:dis}
L_nh(\cos \theta_0)=d\;.
\end{equation}
If $x_0$ is not a Chebyshev node we have
$\theta_0<\theta_{n,k_0}$,\ $0<\sigma_n<1$ and
\begin{equation}\label{eq:ln-nodis}
L_nh(\cos \theta_0)=\sum_{k=1}^{k_0}\ell_{n,k}(\cos \theta_0)\;.
\end{equation}

Let us consider the case where $x_0$ is not a Chebyshev node and observe that
\begin{align*}
\sin((n-1)\theta_0)&=(-1)^{k_0-1}\sin((n-1)\theta_0-(k_0-1)\pi)
\\
&=(-1)^{k_0-1}\sin\left(\pi\left(\frac{n-1}{\pi}\theta_0-(k_0-1)\right)\right)
\\
&=(-1)^{k_0-1}\sin(\pi \sigma_n).
\end{align*}
Then we can rewrite $L_nh(x_0)$ in the following way
\begin{align*}
L_n h(x_0)&=\sum^{k_0}_{k=1}\frac{(-1)^{-k}}{(n-1)(1+\delta_{k,1})}\,\frac{\sin((n-1)\theta_0)\sin \theta_0}
{\cos\theta_0-\cos\theta_{n,k}}
\\
&=
\frac{\sin ((n-1)\theta_0)\sin\theta_0}{2(n-1)(\cos\theta_0-1)}
+\sum^{k_0}_{k=1}\frac{(-1)^{-k}}{(n-1)}\,\frac{\sin((n-1)\theta_0)\sin \theta_0}
{\cos\theta_0-\cos\theta_{n,k}}
\\
&=
\frac{\sin ((n-1)\theta_0)\sin\theta_0}{2(n-1)(\cos\theta_0-1)}
+\frac{\sin(\pi\sigma_n)}
{n-1}\sum^{k_0-1}_{m=0}(-1)^{m}\frac{\sin\theta_0}
{\cos\theta_{n,k_0-m}-\cos\theta_0}
\\
&=\frac{\sin ((n-1)\theta_0)\sin\theta_0}{2(n-1)(\cos\theta_0-1)}
+\frac{\sin(\pi\sigma_n)}{\pi}\sum^{k_0-1}_{m=0}\frac{(-1)^{m}}{\sigma_n+m}
\\
&\quad+\frac{\sin(\pi\sigma_n)}{n-1}\sum^{k_0-1}_{m=0}(-1)^m
\left[
\frac{\sin\theta_0}{\cos\theta_{n,k_0-m}-\cos \theta_0}-
\frac{n-1}{\pi(\sigma_n+m)}\right]
\\
&=\frac{\sin ((n-1)\theta_0)\sin\theta_0}{2(n-1)(\cos\theta_0-1)}
+\frac{\sin(\pi\sigma_n)}{\pi}\sum^{k_0-1}_{m=0}\frac{(-1)^{m}}{\sigma_n+m}
\\
&\quad+\frac{\sin(\pi\sigma_n)}{n-1}\sum^{k_0-1}_{m=0}(-1)^m
\left[
\frac{\sin\theta_0}{\cos\theta_{n,k_0-m}-\cos \theta_0}-
\frac{1}{\theta_0-\theta_{n,k_0-m}}\right],
\end{align*}
where
\[
\theta_0-\theta_{n,k_0-m}=\theta_0-\frac{k_0-m-1}{n-1}\pi=\theta_0-\theta_{n,k_0}+\frac{m}{n-1} \pi
=\frac{\pi}{n-1}(\sigma_n+m).
\]
If we consider the function
\[
g_{\theta_0}(x):=\frac{\sin\theta_0}{\cos x-\cos\theta_0}-\frac{1}{\theta_0-x},\ \ x\in[0,\theta_0[,
\]
we can write
\begin{align}\label{eq:L_nh-no-nodo}
L_nh(x_0)=&\frac{\sin ((n-1)\theta_0)\sin\theta_0}{2(n-1)(\cos\theta_0-1)}
+\frac{\sin(\pi\sigma_n)}{\pi}\sum^{k_0-1}_{m=0}\frac{(-1)^{m}}{\sigma_n+m}
\\
&+\frac{\sin(\pi\sigma_n)}{n-1}\sum^{k_0-1}_{m=0}(-1)^m
 g_{\theta_0}(\theta_{n,k_0-m}). \nonumber
\end{align}
The function $g_{\theta_0}$ is monotone decreasing and bounded since
$
g_{\theta_0}(0)=\frac{\sin\theta_0}{1-\cos\theta_0}-\frac{1}{\theta_0}
$
and
\[
\lim_{x\rightarrow \theta_0^-}g_{\theta_0}(x)=\frac 12 \cot(\theta_0).
\]
For all $n\ge1$ and $\sigma \in [0,1[$, consider the function $f_n:[0,1[\rightarrow\R$ defined by setting
\begin{equation*}
f_{n}(\sigma):=
\left\{
\begin{array}{ll}
\displaystyle \frac{\sin ((n-1)\theta_0)\sin\theta_0}{2(n-1)(\cos\theta_0-1)}
+\frac{\sin(\pi\sigma)}{\pi}\sum^{k_0-1}_{m=0}
  \frac{(-1)^{m}}{\sigma+m}&\\
\qquad\displaystyle
+\frac{\sin(\pi\sigma)}{n-1}\sum^{k_0-1}_{m=0}(-1)^m
g_{\theta_0}(\theta_{n,k_0-m})\;,
&\text{ if }\sigma \in ]0,1[\;,\\
d\;,
&\text{ if }\sigma=0\;;
\end{array}
\right.
\end{equation*}
taking into account \eqref{eq:dis}, \eqref{eq:ln-nodis} and
\eqref{eq:L_nh-no-nodo} we have $L_n h(\cos \theta_0)=f_{n}(\sigma_n)$.

For all $\sigma \in]0,1[$
\begin{align*}
&|f_n(\sigma)-g(\sigma)|
\\
&\quad\leq\left|\frac{\sin(\pi\sigma)}{\pi}\sum^{\infty}_{m=k_0}\frac{(-1)^m}{\sigma+m} \right|
+\frac{1}{2(n-1)(1-\cos \theta_0)}
+\frac{\sin(\pi \sigma)}{n-1}
\left(\frac{1}{\theta_0}+|g_{\theta_0}(\theta_{n,k_0})| \right)
\\
&\quad \leq \frac{\sin(\pi \sigma)}{\pi}\left|\frac{(-1)^{k_0}}{\sigma + k_0}\right|
+\frac{1}{2(n-1)(1-\cos \theta_0)}
+\frac{\sin(\pi \sigma)}{n-1}
\left(\frac{1}{\theta_0}+|g_{\theta_0}(\theta_{n,k_0})| \right)
\\
\\
&\quad \leq \frac{1}{\pi k_0}+
\frac{1}{n-1}
\left(\frac{1}{\theta_0}+|g_{\theta_0}(\theta_{n,k_0})| \right);
\end{align*}
the right-hand side is independent of $\sigma \in ]0,1[$ and
it converges to $0$ as $n\rightarrow \infty$ since
\[
\lim_{n\rightarrow \infty}g_{\theta_0}(\theta_{n,k_0})=
\lim_{x\rightarrow \theta_0^-}g_{\theta_0}(x)=\frac{1}{2}\cot(\theta_0)<\infty\;.
\]
Then we can conclude that the sequence $(f_n)_{n\ge1}$ converges uniformly on $[0,1[$ to the
function $\tilde g:[0,1[\rightarrow\R$ defined as follows
\begin{equation*}
\tilde{g}(x):=\left\{
\begin{array}{ll}
g(x)\;,\qquad&\text{if }x\in ]0,1[\;,
\\
d\;,&\text{if }x=0\;.
\end{array}
\right.
\end{equation*}

Now, we will construct $q$ subsequences
$\left(L_{k_m(n)}h(x_0)\right)_{n\ge1}$, $m=0,\dots,q-1$, of
$(L_{n}h(x_0))_{n\ge1}$ with density $\frac 1q$ such that
\[
\lim_{n\rightarrow\infty}L_{k_m(n)}h(x_0)=
\tilde{g}\left( \frac{m}{q}\right)
 \text{ for all } m=0,\dots, q-1\;.
\]

Fix $m=0,\dots,q-1$; since $GCD(p,q)=1$ we can set $k_m(n):=l+nq+1$, where $l\in \{\,0,\dots,q-1\, \}$ is such that
$lp\equiv m  \mod q$, that is there exists $s\in \Z$ such that
$lp=sq+m$.

So, consider $\left(L_{k_m(n)}h(x_0)\right)_{n\ge1}$ and observe that
for all $m=0,\dots,q-1$, we have $\delta(\{\,k_m(n)\mid n\in \N\, \})=\frac 1q$. It follows, for all $n\ge1$
\begin{align*}
\sigma_{k_m(n)}&=(k_m(n)-1)\frac{p}{q}-
  \left[ (k_m(n)-1)\frac{p}{q}\right]\\&=(l+nq)\frac{p}{q}-
  \left[ (l+nq)\frac{p}{q}\right]=\frac{sq+m+nqp}{q}-\left[\frac{sq+m+nqp}{q}\right]\\
&=s+np+\frac{m}{q}-\left[s+np+\frac{m}{q}\right]=\frac mq
\end{align*}
since $s,np \in \Z$, while $0\leq \frac{m}{q}< 1$.
Then
\begin{align*}
\lim_{n\rightarrow \infty}&L_{k_m(n)}h(x_0)
=\lim_{n\rightarrow \infty}f_n\left( \sigma_{k_m(n)}\right)=\lim_{n\rightarrow \infty}
f_n \left( \frac{m}{q} \right)=\tilde{g}\left( \frac{m}{q} \right)\;.
\end{align*}

Therefore, by \cite[Proposition 1.6]{CMT}
we have that
for all $m=1,\dots,q$
\[
i\left(L_nh(x_{0}),
\tilde{g}\left(\frac{m}{q}\right)\right)
\geq \frac 1q\;.
\]
Now, we have $q$ different statistical limits with index $\frac 1q$, so by  \cite[Proposition 1.7]{CMT}
it necessarily follows
\[
i\left(L_nh(x_{0});
\tilde{g}\left(\frac{m}{q}\right)\right)
=\frac 1q.
\]
This completes the proof of part i).

The case where $\frac{\theta _0}{\pi}$ is irrational is similar to the analogous case considered in the proof of
\cite[Theorem 2.1 ii)]{CMT}
\qed

We observe that if $g\left(\frac{m_0}{q}\right)=d$ for some $m_0=1,\dots,q-1$ then item i) becomes
  \[
    i\left( L_nh(x_0);d\right)=\frac2q\;,\qquad
    i\left(L_nh(x_0);
    g\left( \frac{m}{q}\right)\right)=\frac1q
  \]
for every $m=1,\dots, q-1,\ m\ne m_0$.

At this point, we extend Theorem \ref{th:Lagrange} to a larger classes of
functions, namely on the space
$\mathcal{C}+H$ where $\mathcal{C}$ denotes the space of all $f\in C([-1,1])$ such that $f$ is either  monotone
 on $[-1,1]$ or $f$ satisfies the Dini-Lipschitz condition $\omega(f,\delta)=o(|\log\delta|^{-1})$,
and $H$ is the linear space generated by
\[
\{\,h_{x_0,d}\mid x_0\in]-1,1[,\ d\in\R\, \}\;.
\]

Observe that if $f\in \mathcal{C}+H$ there exists at most a finite number of points $x_1,\dots,x_N$
of discontinuity with finite left and right limits $f(x_i-0)$ and $f(x_i+0)$,\ $i=1,\dots,N$.

Then we can state the following theorem.

\begin{theorem} \label{th:Lagrange_extended}
Let $f\in \mathcal{C}+H$ with a finite number $N$ of
points of discontinuity of the first kind at $x_1,\dots,x_N\in]-1,1[$.
For every $i=1,\dots,N$ consider $\theta_i\in]0,\pi[$ such that $x_i=\cos
\theta_i$,\ $d_i:=f(x_i)$ and define the function
\[
g_i(x):=f(x_i-0)+(f(x_i+0)-f(x_i-0))g(x)\;.
\]

Then, the sequence $\left( L_n f\right)_{n\ge1}$ converges uniformly
to $f$
on every compact subset of $]-1,1[\setminus\{\,x_1,\dots,x_N\, \}$.

Moreover for all $i=1,\dots,N$ the sequence $(L_n f(x_i))_{n\ge1}$ has the following
behavior
\begin{itemize}
\item[i)] if $\frac{\theta_i}{\pi}=\frac pq $ with $p,q\in \N$, $q\neq 0$ and
$GCD(p,q)=1$, then \[
    i\left(L_nh(x_i);d_i\right)=\frac1q\;,\qquad
    i\left(L_nh(x_i);
    g_i\left( \frac{m}{q}\right)\right)=\frac1q\;,
    \quad m=1,\dots, q-1.
  \]

\item [ii)] if $\frac{\theta_i}{\pi}$ is irrational and if
$A\subset \R$ is a Peano-Jordan measurable set,
then
\[
i\left(L_nh(x_i);A\right)=|g_i^{-1}(A)|\;,
\]
where $|\cdot|$ denotes the Peano-Jordan measure.
\end{itemize}
\end{theorem}

\proof
We assume $x_1<\dots<x_N$. We can write $f=F+\sum _{k=1}^{N}c_kh_k$,
where $F\in \mathcal{C}$ and
$h_i:=h_{x_i,\tilde d_i}$ for every $i=1,\dots,N$.

Since $F$ is continuous we have
\begin{align*}
f(x_i+0)-\sum_{k=1}^{i-1}c_k-c_i=F(x_i+0)=F(x_i-0)=f(x_i-0)-\sum_{k=1}^{i-1}c_k\;,
\end{align*}
from which
\[
c_i=f(x_i+0)-f(x_i-0)
\]
and
\begin{equation}\label{eq:F(x_i)}
 F(x_i)=f(x_i-0)-\sum_{k=1}^{i-1}c_k\;.
\end{equation}
Moreover
\begin{align*}
d_i=f(x_i)&=F(x_i)+\sum^{i-1}_{k=1}c_k h_k(x_i)+c_i \tilde d_i
\\
&=F(x_i)+\sum^{i-1}_{k=1}c_k +\left(f(x_i+0)-f(x_i-0)\right)\tilde d_i
\\
&=f(x_i-0)+\left(f(x_i+0)-f(x_i-0)\right)\tilde d_i\;.
\end{align*}
and hence
\[
\tilde d_i=\frac{d_i-f(x_i-0)}{f(x_i+0)-f(x_i-0)}\;.
\]

The first part of our statement is a trivial consequence of the linearity of
Lagrange interpolation operators. Indeed
$F\in \mathcal{C}$
and therefore $L_n F\to F$ uniformly in compact subsets of $]-1,1[$
(see e.g. \cite[Theorem 3.2, p. 24]{V} and \cite{MM} in the case in which $F$ is monotone,
while we refer to
\cite[Theorem 14.4, p. 335]{S} in the case $F$ satisfies the Dini-Lipschitz condition);
moreover for every $k=1,\dots,N$, by Theorem \ref{th:Lagrange}
$L_n h_k\to h_k$ converges uniformly to $h_k$ on compact subsets
of $[-1,1]\setminus \{\,x_k\, \}$. Then $L_n f=L_n F+\sum^{N}_{k=1}c_k L_n h_k$ converges
uniformly to $f$ on compact subsets of $]-1,1[\setminus\{\,x_1,\dots,x_N\, \}$.

Now we establish property i). We fix a point $x_i$ of discontinuity
and following the same line of the proof of Theorem \ref{th:Lagrange}
we construct the subsequences $(k_m(n))_{n\ge1}$,\ $m=1,\dots,q$.
Since
\[
L_{k_m(n)}f(x_i)=L_{k_m(n)} F(x_i)+
\sum_{
\genfrac{}{}{0pt}{}{k=1}{k\neq i}
}^{N}c_k L_{k_m(n)} h_k(x_i)+c_i
L_{k_m(n)} h_i(x_i)
\]
and taking into account \eqref{eq:F(x_i)} and that $F\in \mathcal{C}$,
from Theorem \ref{th:Lagrange}
the right-hand side converges to
\begin{align*}
&F(x_i)+\sum^{i-1}_{k=1}c_k h_k(x_i)+c_i g_i\left(\frac{m}{q}\right) \\
&\qquad =f(x_i-0)+\left(f(x_i+0)-f(x_i-0)\right)g_i\left(\frac{m}{q}\right)\\
&\qquad =g_i\left(\frac{m}{q}\right)
\end{align*}
for $m=0,\dots,q-1$.

Finally, we prove property ii).
For every $i=1,\dots,N$ we have
\[
L_{n}f(x_i)=L_{n} F(x_i)+
\sum_{
\genfrac{}{}{0pt}{}{k=1}{k\neq i}
}^{N}
c_k L_{n} h_k(x_i)+c_i
L_{n} h_i(x_i)\;.
\]
For the sake of simplicity let us denote
\[
y_n:=L_n f(x_i)\;,\quad
z_n:=L_n F(x_i)+
  \sum_{\genfrac{}{}{0pt}{}{k=1}{k\neq i}}^{N}c_k L_n h_k(x_i)\;,\quad
x_n:=c_iL_n h_i(x_i)
\]
(thus $y_n=z_n+x_n$) and
\[z:=F(x_i)+\sum_{k=1}^{i-1}c_k
h_k(x_i)=f(x_i-0)
\]
(see \eqref{eq:F(x_i)}).

Since \eqref{eq:F(x_i)} and that $F\in \mathcal{C}$
we can apply \cite[Theorem 3.2, p. 24]{V}
(or \cite[Theorem 14.4, p. 335]{S}) and from
Theorem \ref{th:Lagrange} we obtain $z_n\rightarrow z$ and
moreover $i(c_i^{-1}x_n;A)=|g^{-1}(A)|$ for every bounded Peano-Jordan measurable
set $A\subset \R$. Hence $i(x_n;A)=|g^{-1}(c_i^{-1}A)|$, that is

\[|g^{-1}(c_i^{-1}A)|
=\inf_{\varepsilon>0}
\delta_-(\{\,n\in \N\mid x_n\in A+B_\varepsilon\, \})\;.\]

Fix $\varepsilon>0$; if $x_n\in A+B_\varepsilon$, from the equality $x_n=y_n-z_n$ we get
\[
y_n\in A+B_\varepsilon+z_n=A+B_\varepsilon+z+z_n-z\;.
\]
Now, let $\nu \in \N$ such that $|z_n-z| < \varepsilon$ for all $n\geq \nu$,
then for every
$n\geq \nu$ we have $z_n-z\in B_\varepsilon$ and consequently $y_n\in A+B_{2\varepsilon}+z$.
Therefore
\begin{align*}
\{\,n\geq \nu\mid x_n\in A+B_\varepsilon\, \}\subset \{\,n\geq \nu\mid y_n\in
A+B_{2\varepsilon}+z\, \},
\end{align*}
that is
\begin{align}\label{eq:111}
\delta_- (\{\,n\in \N\mid x_n\in A+B_\varepsilon\, \})
\leq
\delta_-( \{\,n\in \N\mid y_n\in A+B_{2\varepsilon}+z\, \})\;.
\end{align}

On the other hand, if $y_n\in A+B_{2\varepsilon}+z$,
then $x_n=y_n-z_n\in A+B_{2\varepsilon}+z-z_n$.
In this case for every $n\geq \nu$, we have
$z-z_n\in B_\varepsilon$ and therefore
$x_n\in A+B_{3\varepsilon}$; hence
\begin{align}\label{eq:222}
\delta_-
(\{\,n\in \N\mid x_n\in A+B_{3\varepsilon}\, \})\geq
\delta_- \{\,n\in \N\mid y_n\in A+B_{2\varepsilon}+z\, \}\;.
\end{align}

Taking the infimum over $\varepsilon>0$
in \eqref{eq:111} and \eqref{eq:222}
we can conclude that
$i(x_n,A)\leq i(y_n,A+z)\leq i(x_n,A)$ which yields
\[
i(y_n,A+z)=i(x_n,A)=|g^{-1}(c_i^{-1}A)|\;.
\]
We conclude that
$i(y_n,A)=|g^{-1}(c_i^{-1}(A-z))|
=\left|g^{-1}\left(\frac{A-f(x_i-0)}{f(x_i+0)-f(x_i-0)}\right)\right|=|g_i^{-1}(A)|$
for every Peano-Jordan measurable set $A\subset \R$.
\qed

As in Theorem \ref{th:Lagrange} item i) becomes
  \[
    i\left( L_nh(x_0);d_i\right)=\frac2q\;,\qquad
    i\left(L_nh(x_0);
    g_i\left( \frac{m}{q}\right)\right)=\frac1q
  \]
for every $m=1,\dots, q-1,\ m\ne m_0$, if $g_i\left(\frac{m_0}{q}\right)=d_i$ for some $m_0=1,\dots,q-1.$


Now, we can consider the bivariate Lagrange interpolation polynomials
$(L_{n,m})_{n,m\geq 1}$ on the Chebyshev nodes of second kind plus the endpoints $\pm 1$ defined by
\begin{equation}\label{eq:Lagrange-bivariate}
L_{n,m}(f)(x,y)=\sum_{i=1}^n\sum_{j=1}^m
    \omega _{n,i}(x)\omega_{m,j}(j)f(x_{n,i},y_{m,j})
\end{equation}
where $f$ is a suitable function defined on $[-1,1]^2$ and
\[
x_{n,i}=\cos\frac{i-1}{n-1}\pi\;,\quad i=1,\dots,n\;.
\]
Moreover, setting $x=\cos \theta$,
\[
\omega _{n,i}(x)=\frac{(-1)^i}{(1+\delta_{i,1}+\delta_{i,n}(n-1)}
\frac{\sin \theta \sin ((n-1)\theta)}{x-x_{n,i}}\;.
\]
Consider $z_0=(x_0,y_0)\in]-1,1[\times]-1,1[$ and define the following function $h_{z_0}:\, [0,1]\times [0,1]\rightarrow \R$,
\begin{equation} \label{eq:def_hS1}
h_{z_0}(x,y):=\left\{
\begin{array}{ll}
1\;,\quad& (x,y)\in [x_0,1]\times [y_0,1]\;,\\
0\;,\quad& (x,y)\in [-1,1]\times[-1,1]\setminus[x_0,1]\times [y_0,1]\;.
\end{array} \right.
\end{equation}
In order to state the convergence properties of the sequence $(L_{n,m,s}h_{z_0})_{n\ge1}$,
 we consider the function $G:\, ]0,1[\times ]0,1[\rightarrow \R$ defined 
\[
G(x,y):=g(x)g(y)\;,
\]
where $g:]0,1[\rightarrow\R$ is the function defined in \eqref{eq:g}.

In the following result we describe the behavior of the bivariate Lagrange polynomials evaluated at the function
$h_{z_0}$, using the index of convergence for double sequences of real numbers.

\begin{theorem} \label{th:lagrange-multi}
Let $z_0=(x_0,y_0)=(\cos \theta_0,\cos \gamma_0)\in]-1,1[\times]-1,1[$ and $h:=h_{z_0}$ be defined by \eqref{eq:def_hS1}.
Then the sequence $( {L}_{n,m,s}h)_{n\ge1}$ converges uniformly to $h$
on every compact subset of $[-1,1]\times [-1,1]\setminus Q$, where $Q:=([x_0,1]\times \{\,y_0\, \})
\cup (\{\,x_0\, \}\times [y_0,1])$.

As regards the behaviour of the sequence $({L}_{n,m,s}h(x,y))_{n,m\ge1}$
where $(x,y)\in Q$, we have:
\begin{itemize}

\item[1)] if $x=x_0$ and $y\in]y_0,1]$, we have to consider the following cases:
\begin{itemize}
\item[i)] if $\frac{\theta_0}{\pi}=\frac{p_1}{q_1}$ with $p_1,q_1\in \N$, $q_1\neq 0$,
$GCD(p_1,q_1)=1$, then
\begin{align*}
i\left({L}_{n,m}h(x_0,y);g\left(\frac{m_1}{q_1}\right)\right)=
\frac{1}{q_1}
\;,\; m_1=0,\dots,q_1-1\;;
\end{align*}
\item[ii)] if $\frac{\theta_0}{\pi}$ is irrational, we have
\[
i\left({L}_{n,m}h(x_0,y);A\right)=|g^{-1}(A)|
\]
for every Peano-Jordan measurable set $A\subset \R$.
\end{itemize}

\item[2)] If $x\in]x_0,1]$ and $y=y_0$, we have to consider the following cases:
\begin{itemize}
\item[i)] if $\frac{\gamma_0}{\pi}=\frac{p_2}{q_2}$ with $p_2,q_2\in \N$, $q_2\neq 0$,
$GCD(p_2,q_2)=1$, then
\begin{align*}
i\left({L}_{n,m}h(x,y_0);g\left(\frac{m_2}{q_2}\right)\right)=
\frac{1}{q_2}
\;,\; m_2=0,\dots,q_2-1\;;
\end{align*}

\item[ii)] if $\frac{\gamma_0}{\pi} $ is irrational, we have
\[
i\left({L}_{n,m}h(x,y_0);A\right)=|g^{-1}(A)|
\]
for every Peano-Jordan measurable set $A\subset \R$ if $s>1$\;.
\end{itemize}

\item[3)] If $x=x_0$ and $y=y_0$, we have to consider the following cases:
\begin{itemize}
\item[i)] if $\frac{\theta_0}{\pi}=\frac{p_1}{q_1}$ and $\frac{\gamma_0}{\pi}=\frac{p_2}{q_2}$ with $p_i,q_i\in \N$, $q_i\neq 0$,
$GCD(p_i,q_i)=1$,\ $i=1,2$, then
\begin{align*}
i\left( {L}_{n,m}h(x_0,y_0);
g\left(\frac{m_1}{q_1}\right)g\left(\frac{m_2}{q_2}\right)\right)=
\frac{1}{q_1q_2}
\end{align*}
for $m_1=0,\dots,q_1-1$ and $m_2=0,\dots,q_2-1$;

\item[ii)] if $\frac{\theta_0}{\pi}=\frac{p_1}{q_1}$ with $p_1,q_1\in \N$, $q_1\neq 0$,
$GCD(p_1,q_1)=1$ and $\frac{\gamma_0}{\pi}$ is irrational, then
\begin{align*}
&i\left({L}_{n,m}h(x_0,y_0);
\left[0,g\left(\frac {j}{q_1}\right)\right]\right)\geq\frac{1}{q_1},\qquad j=0,\dots,q_1-1\;;
\end{align*}
\item[iii)] if $\frac{\theta_0}{\pi}$ is irrational and $\frac{\gamma_0}{\pi}=\frac{p_2}{q_2}$ with $p_2,q_2\in \N$, $q_2\neq 0$,
$GCD(p_2,q_2)=1$, then
\begin{align*}
&i\left({L}_{n,m}h(x_0,y_0);\left[0,g\left(\frac{j}{q_2}\right)\right]\right)
\geq\frac{1}{q_2},\qquad j=0,\dots,q_2-1\;;
\end{align*}

\item[iv)] if $\frac{\theta_0}{\pi}$ and $\frac{\gamma_0}{\pi}$ are both irrational, then
\[
i\left({L}_{n,m}h(x_0,y_0);A\right)=|G^{-1}(A)|
\]
for every Peano-Jordan measurable set $A\subset \R$.
\end{itemize}
\end{itemize}

\end{theorem}
\proof
We define the functions $h_1,\;h_2:[-1,1]\to \R$ as follows
\begin{equation}\label{eq:h1-h2}
h_1(t):=\left\{
\begin{array}{ll}
0\, &\text{if }t< x_0,
\\
1\, &\text{if }t\geq x_0,
\end{array}
\right.
\qquad
h_2(t):=\left\{
\begin{array}{ll}
0\, &\text{if }t< y_0,
\\
1\, &\text{if }t\geq y_0\;,
\end{array}
\right.
\end{equation}
then we can write
\[
h(x,y)=h_1(x)h_2(y)\;,
\]
and consequently we have
\begin{align}\label{eq:decomp-Lnm}
&L_{n,m}h(x,y)=L_nh_1(x)L_mh_2(y)\;.
\end{align}
Moreover, by Theorem \ref{th:Lagrange},
it follows that
\begin{equation}\label{eq:conv_compact1}
\lim_{n\rightarrow \infty}L_nh_1=h_1 \text{ uniformly on }[-1,1]\setminus\{\,x_0\, \}
\end{equation}
and
\begin{equation}\label{eq:conv_compact2}
\lim_{m\rightarrow \infty}L_mh_2=h_2 \text{ uniformly on }[-1,1]\setminus\{\,y_0\, \}\;.
\end{equation}

Let us consider a compact set $K\subset [-1,1]\times[-1,1]\setminus Q$; there exist
$-1\leq a_1<x_0< a_2\leq 1$ and $-1\leq b_1<y_0< b_2\leq1$ such that
\[
K\subset [-1,1]\times[-1,b_1]\cup [-1,a_1]\times [-1,1]\cup [a_2,1]\times[b_2,1].
\]
First let us consider $(x,y)\in [-1,1]\times[-1,b_1]$ and write $x=\cos \theta$ and $y=\cos \gamma$;
then
\begin{equation}\label{eq:|L_n,mh|}
|L_{n,m}h(x,y)|\leq |L_nh_1(x)||L_mh_2(y)|.
\end{equation}
For sufficiently large $n,m\geq 1$ there exist $k_0,\ell_0\geq 1$ such that
\[
x_{n,k_0+1}<x_0\leq x_{n,k_0},\qquad y_{m,\ell_0+1}<y_0\leq y_{m,\ell_0}\;,
\]
then we can observe that
\begin{align}\label{eq:estimate1_L_nh_1}
|L_nh_1(x)|&\leq \frac{|\sin (n-1)\theta \sin \theta|}{n-1}
\left|\sum^{k_0}_{k=1}\frac{(-1)^k}{(1+\delta_{k,1})(\cos \theta-\cos \theta_{n,k})}\right|
\nonumber
\\
&\leq \frac{|\sin (n-1)\theta \sin \theta|}{2(n-1)|\cos \theta - \cos \theta_{n,1}|}
+\frac{|\sin (n-1)\theta \sin \theta|}{(n-1)|\cos \theta - \cos \theta_{n,k_0}|}
\nonumber
\\
&=\frac{|\omega_{n,1}(x)|}{2}+|\omega_{n,k_0}(x)|
\end{align}
where, in the last inequality, we have used the fact that the following function
\[
t\in [0,\pi]\rightarrow \frac{1}{\cos \theta -\cos t}
\]
is monotone decreasing.
\\
Let us observe that for $k=1,\dots,k_0$ we have
\begin{align}\label{eq:estimate_omega}
|\omega_{n,k}(x)|&=\left|\frac{\sin (n-1)\theta \sin \theta}{(n-1)(\cos \theta - \cos \theta_{n,k})}\right|
=\left|\frac{\sin\left((n-1)(\theta-\theta_{n,k})+(k-1)\pi\right)\sin \theta}
{2(n-1)\sin \left(\frac{\theta -\theta_{n,k}}{2}\right)\sin \left(\frac{\theta+\theta_{n,k}}{2}\right)}\right|
\nonumber
\\
&=\left|\frac{\sin\left((n-1)(\theta-\theta_{n,k})\right)\sin \theta}
{2(n-1)\sin \left(\frac{\theta -\theta_{n,k}}{2}\right)\sin \left(\frac{\theta+\theta_{n,k}}{2}\right)}\right|
\nonumber
\\
&=\frac{2\sin\theta}
{\theta+\theta_{n,k}}\frac{|\sin\left((n-1)(\theta-\theta_{n,k})\right)|}{(n-1)|\theta-\theta_{n,k}|}
\frac{\frac{|\theta-\theta_{n,k}|}{2}}{\left|\sin \left(\frac{\theta -\theta_{n,k}}{2}\right)\right|}
\frac{\frac{\theta+\theta_{n,k}}{2}}{\left|\sin \left(\frac{\theta+\theta_{n,k}}{2}\right)\right|}
\nonumber
\\
&\leq \frac{2\sin\theta}
{\theta}\frac{|\sin\left((n-1)(\theta-\theta_{n,k}|\right)|}{(n-1)|\theta-\theta_{n,k}|}
\frac{\frac{|\theta-\theta_{n,k}|}{2}}{\sin \left|\frac{\theta -\theta_{n,k}}{2}\right|}
\frac{\frac{\theta+\theta_{n,k}}{2}}{\left|\sin \left(\frac{\theta+\theta_{n,k}}{2}\right)\right|}
\nonumber
\\
&\leq 2C
\end{align}
where the existence of the constant $C\geq 1$ is a consequence of the boundedness of the functions
$\sin\alpha/\alpha$ on $[0,n\pi]$ for all $n\geq 1$ and $\alpha/\sin\alpha$ on
$[a,b]\subset [0,\pi[$;
in particular, observe that $\frac{\sin\alpha}{\alpha}\leq 1$ while
$\frac{\alpha}{\sin\alpha}\geq 1$ and
$\left|\frac{\theta-\theta_{n,k}}{2}\right|,\left|\frac{\theta+\theta_{n, k}}{2}\right|\neq \pi$,
since $0<\theta_{n,k}<\pi$ for all $k$; notice also that \eqref{eq:estimate_omega}
does not depend on the particular
choice of $k\geq1$ and it holds for all $x\in [-1,1]$. Then \eqref{eq:estimate1_L_nh_1} becomes
\begin{equation}\label{eq:estimate2_L_nh_1}
|L_nh_1(x)|\leq 3C
\end{equation}
and this estimate is uniform with respect to $n\geq 1$ and $x\in [-1,1]$.

We can conclude that
\[
|L_{n,m}h(x,y)|\leq 3C |L_{m}h_2(y)|
\]
where, by \eqref{eq:conv_compact2}, the last term converges to $0$ as $m\rightarrow \infty$.
So we can conclude that
\[
\lim_{n,m\rightarrow \infty }L_{n,m}h=h\qquad \text{ uniformly in }[-1,1]\times[-1,b_1]\;.
\]
Arguing in a similar way, we can get the uniform convergence of $(L_{n,m}h)_{n,m\geq 1}$ to
$h$ in $[-1,a_1]\times[-1,1]$.

If $(x,y)\in [a_2,1]\times[b_2,1]$, then
\begin{align*}
|L_{n,m}h(x,y)-1|&= |L_nh_1(x)L_mh_2(y)-1|
\end{align*}
which is uniformly convergent to 0 as $n,m\rightarrow \infty$. Therefore we
can conclude that
\[
\lim_{n,m\rightarrow \infty }L_{n,m}h=h\qquad \text{ uniformly in }[a_2,1]\times[b_2,1].
\]

We start with the proof of property 1). Let $y\in ]y_0,1]$.
From equation \eqref{eq:decomp-Lnm} we have
\begin{align*}
&L_{n,m}h(x_0,y)=L_nh_1(x_0)L_mh_2(y)\nonumber\;,
\end{align*}
and thank to \eqref{eq:conv_compact2} we have
\[
\lim_{m\to \infty }
L_{n,m}h(x,y_0)= L_{n}h_1(x_0)\;.
\]
From  \eqref{eq:estimate2_L_nh_1}
we have that $L_nh_1(x_0)$ is bounded, then the previous limit is uniform with respect $n\in\N$.
Therefore we can apply Proposition \ref{pr:1-bivariate}
with $k(m)=m$ (and consequently $\alpha=1$) and Theorem \ref{th:lagrange-multi} with $h$ replaced by $h_1$, and conclude the proof of 1).

The proof of property 2) is at all similar to that of property 1) interchanging the role of $x$ and $y$.

Now we prove 3).
From \eqref{eq:decomp-Lnm} we have
\begin{align*}
&L_{n,m}h(x_0,y_0)=L_nh_1(x_0)L_mh_2(y_0)\;.
\end{align*}

Arguing as in Theorem \ref{th:lagrange-multi} and taking into account that the value $d$ is set to $1$ and $g(0)=1$, we can consider
$\left(L_{r_{i}(n)}h_1(x_0)\right)_{n\geq 1}$,
$i=0,\dots,q_1-1$ and $\left(L_{s_{j}(m)}h_2(y_0)\right)_{m\geq1}$,
$j=0,\dots,q_2-1$, subsequences respectively of $(L_{n}h_1(x_0))_{n\geq1}$ and
$(L_{m}h_2(y_0))_{m\geq1}$ with density respectively $1/q_1$ and $1/q_2$ such
that
\begin{align*}
&\lim_{n\rightarrow \infty}L_{r_{ji}(n)}h_1(x_0)
  =g\left(\frac{i}{q_1}\right),
\quad i=0,\dots,q_1-1,
\end{align*}
and
\begin{align*}
&\lim_{n\rightarrow \infty}L_{s_{j}(m)}h_1(x_0)
  =g\left(\frac{j}{q_2}\right),
\quad j=0,\dots,q_2-1\;.
\end{align*}
Therefore we can consider
$q_1q_2$ subsequences of $\left({L}_{n,m}h(x_0,y_0)\right)_{n,m\geq 1}$,
let us say $\left(L_{r_{i}(n),s_{j}(m)}h(x_0,y_0)\right)_{n,m\geq 1}$ with
$i=0,\dots,q_1-1$ and $j=0,\dots,q_2-1$ such that
$\delta(\{\,r_{i}(n),s_{j}(m)|n,m \in \N \, \})=\frac{1}{q1q_2}$ and
\[
\lim_{n,m\rightarrow \infty}L_{r_{i}(n),s_{j}(m),s}h(x_0,y_0)=G\left(
\frac{i}{q_1},\frac{j}{q_2}\right)
\]
where $i=0,\dots,q_1-1$, $j=0,\dots,q_2-1$. From Propositions \ref{prop:char_index_multi}  and \ref{pr:sum_gen_multi}    we have the result.

Now, let us prove 3) case ii).
Suppose that $\frac{\theta_0}{\pi}=\frac{p_1}{q_1}$, $p_1$, $q_1\in \N$, $q_1\neq 0$,
$GCD(p_1,q_1)=1$ and $\frac{\gamma_0}{\pi}$ is irrational,
from \eqref{eq:decomp-Lnm} we have
\begin{align*}
&L_{n,m}h(x_0,y_0)=L_nh_1(x_0)L_mh_2(y_0)\;.
\end{align*}
We can consider
$\left(L_{r_{i}(n)}h_1(x_0)\right)_{n\geq 1}$,
$i=0,\dots,q_1-1$ subsequences of $(L_{n}h_1(x_0))_{n\geq1}$  with density $1/q_1$ such that
\begin{align*}
&\lim_{n\rightarrow \infty}L_{r_{j}(n),m}h(x_0,y_0)
  =g\left(\frac{i}{q_1}\right)L_{m}h_2(y_0),
\quad i=0,\dots,q_1-1.
\end{align*}
Applying Proposition \ref{pr:1-bivariate} we have
\begin{equation}\label{eq:dis1}
i\left(L_{n,m}h(x_0,y_0);A\right)\geq \frac{1}{q_1}
i\left( d_2L_{m}h_2(y_0);A\right)=\frac{|(d_2g)^{-1}(A)|}{q_1}.
\end{equation}
and
\begin{align}\label{eq:dis2}
&i\left(L_{n,m}h(x_0,y_0);A\right)\geq\frac{1}{q_1}
  i\left(g\left(\frac{i}{q_1}\right)L_{m}h_2(y_0);A\right)\nonumber\\
&\qquad =\frac{|G^{-1}_{i}(A)|}{q_1},\qquad i=1,\dots,q_1-1.
\end{align}
where $G_{i}(t)=G\left(\frac{i}{q_1},t\right)=
g\left(\frac{i}{q_1}\right)g(t)$.
Since the sum of indices can't exceed 1 in
inequalities  \eqref{eq:dis1} and
\eqref{eq:dis2} we have equalities.

The proof of 3) case iii) is at all similar to the previous one interchanging
the role of $x$ and $y$.

Let us conclude the proof of our theorem, considering the case in which both $\frac{\theta_0}{\pi}$ and $\frac{\gamma_0}{\pi}$
are irrational. In this case, we have
$i(L_nh_1(x_0);A)=|g^{-1}(A)|$ and
$i(L_mh_2(y_0);A)=|g^{-1}(A)|$ for every $A$ Peano-Jordan measurable set,
apply Theorem \ref{pr:prod-succ}
to the sequences $(L_{n}h_1(x))_{n\geq 1}$ and to $(L_{m}h_2(y_0))_{m\geq 1}$ and taking into account that
$L_{n,m}h(x_0,y_0)=L_nh_1(x_0)L_mh_2(y_0)$
the claim easily follows.
\qed

\section{Bivariate Shepard operators on discontinuous functions}\label{sc:Shepard}

Among all different kinds of bivariate Shepard operators (see e.g. \cite{G}), for the sake of simplicity
we concentrate our attention to the bivariate Shepard operators obtained as tensor product of univariate Shepard operators
\begin{equation}\label{eq:Shepard_tensor}
S_{n,m,s}f(x,y):=
\sum^n_{i=0}\sum^m_{j=0}\frac{|x-x_i|^{-s}}{\sum^n_{k=0}|x-x_k|^{-s}}\,
\frac{|y-y_j|^{-s}}{\sum^m_{k=0}|y-y_k|^{-s}}f(x_i,y_j),
\end{equation}
where $f$ is a suitable function defined in $[0,1]\times[0,1]$, $s\geq 1$, $n,\ m\geq1$
and $\left((x_i,y_j)\right)_{i,j}$
is the matrix $(n+1)\times(m+1)$ of equispaced nodes in $[0,1]\times[0,1]$, that is
\[
x_i:=\frac in,\ i=0,\dots,n,\qquad y_j:=\frac jm,\ j=0,\dots,m.
\]

The aim of this section is the study of their behavior on
a particular class of bivariate functions having suitable discontinuities defined as follows.

Consider $z_0=(x_0,y_0)\in[0,1]\times[0,1]$ and define the following function $h_{z_0,d}:\, [0,1]\times [0,1]\rightarrow \R$,
\begin{equation} \label{eq:def_hS}
h_{z_0,d}(x,y):=\left\{
\begin{array}{ll}
1\;,\quad& (x,y)\in [0,x_0]\times [0,y_0]\;,\\
0\;,\quad& \text{otherwise}\;,
\end{array} \right.
\end{equation}
In order to state the convergence properties of the sequence $(S_{n,m,s}h_{z_0,d})_{n\ge1}$,
for every $s>1$ we consider the function $g_s:[0,1[\rightarrow\R$ defined as follows
\[
g_s(t):=\left\{
\begin{array}{ll}
\displaystyle\frac{\zeta(s,t)}{\zeta(s,t)+\zeta(s,1-t)}\;,\qquad &t\in]0,1[\;,\\
1\;,&t=0\;,
\end{array} \right.
\]
where $\zeta$ denotes the Hurwitz zeta function:
\begin{equation} \label{eq:def_H}
\zeta(s,a):=\sum^{+\infty}_{n=0}\frac{1}{(n+a)^s}
\end{equation}
for all $s,a\in \C$ such that $\Re[s]>1$ and $\Re[a]>0$. The previous series is absolutely convergent and
its sum can be extended to a meromorphic function defined for all $s\neq 1$.

We consider also $G_s:\, [0,1[\times [0,1[\rightarrow \R$ defined as follow
\[
G_s(x,y):=g_s(x)g_s(y).
\]

We have the following result.

\begin{theorem} \label{th:shepard}
Let $z_0=(x_0,y_0)\in[0,1]\times[0,1]$ and $h:=h_{z_0,d}$ be defined by \eqref{eq:def_hS}.
Then for every $s\ge1$ the sequence $\left( S_{n,m,s}h\right)_{n\ge1}$ converges uniformly to $h$
on every compact subset of $[0,1]\times [0,1]\setminus Q$, where $Q:=([0,x_0]\times \{\,y_0\, \})
\cup (\{\,x_0\, \}\times [0,y_0])$ .

As regards the behavior of the sequence $(S_{n,m,s}h(x,y))_{n,m\ge1}$
where $(x,y)\in Q$, we have:
\begin{itemize}
\item[1)] if $x\in[0,x_0[$ and $y=y_0$, we have to consider the following cases:
\begin{itemize}
\item[i)] if $y_0=\frac{p_2}{q_2}$ with $p_2,q_2\in \N$, $q_2\neq 0$,
$GCD(p_2,q_2)=1$, then
\begin{align*}
i\left(S_{n,m,s}h(x,y_0);g_s\left(\frac{m_2}{q_2}\right)\right)=
\frac{1}{q_2}
\;,\; m_2=0,\dots,q_2-1\;
\end{align*}
if $s>1$; while
\begin{align*}
i\left(S_{n,m,s}h(x,y_0);1\right)&=\frac{1}{q_2}\;,\quad
i\left(S_{n,m,s}h(x,y_0);\frac12\right)=
1-\frac{1}{q_2}\;
\end{align*}
if $s=1$;
\item[ii)] if $y_0$ is irrational, we have
\[
i\left(S_{n,m,s}h(x,y_0);A\right)=|g_s^{-1}(A)|
\]
for every Peano-Jordan measurable set $A\subset \R$ if $s>1$; while
\[
i\left(S_{n,m,s}h(x,y_0);\frac12\right)=1
\]
if $s=1$.
\end{itemize}

\item[2)] If $x=x_0$ and $y\in[0,y_0[$, we have to consider the following cases:
\begin{itemize}
\item[i)] if $x_0=\frac{p_1}{q_1}$ with $p_1,q_1\in \N$, $q_1\neq 0$,
$GCD(p_1,q_1)=1$, then
\begin{align*}
i\left(S_{n,m,s}h(x_0,y);g_s\left(\frac{m_1}{q_1}\right)\right)=
\frac{1}{q_1}
\;,\; m_1=0,\dots,q_1-1\;,
\end{align*}
if $s>1$; while
\begin{align*}
i\left(S_{n,m,s}h(x_0,y);1\right)&=\frac{1}{q_1}\;,\quad
i\left(S_{n,m,s}h(x_0,y);\frac12\right)=
1-\frac{1}{q_1}\;
\end{align*}
if $s=1$;
\item[ii)] if $x_0$ is irrational, we have
\[
i\left(S_{n,m,s}h(x_0,y);A\right)=|g_s^{-1}(A)|
\]
for every Peano-Jordan measurable set $A\subset \R$ if $s>1$; while
\[
i\left(S_{n,m,s}h(x_0,y);\frac12\right)=1
\]
if $s=1$.
\end{itemize}

\item[3)] If $x=x_0$ and $y=y_0$, we have to consider the following cases:
\begin{itemize}
\item[i)] if $x_0=\frac{p_1}{q_1}$ and $y_0=\frac{p_2}{q_2}$ with $p_i,q_i\in \N$, $q_i\neq 0$,
$GCD(p_i,q_i)=1$,\ $i=1,2$, then
\[
i\left(S_{n,m,s}h(x_0,y_0);G_s\left( \frac{m_1}{q_1},\frac{m_2}{q_2}\right)\right)=
\frac{1}{q_1q_2},
\]
where $m_1=0,\dots, q_1-1$, $m_2=0,\dots q_2-1$, if $s>1$.
While, if $s=1$,
\begin{align*}
i\left(S_{n,m,s}h(x_0,y_0);\frac12\right)&=\frac{1}{q_1q_2}\;,\quad
i\left(S_{n,m,s}h(x_0,y_0);\frac14\right)=
1-\frac{1}{q_1q_2}\;;
\end{align*}

\item[ii)] if $x_0=\frac{p_1}{q_1}$ with $p_1,q_1\in \N$, $q_1\neq 0$,
$GCD(p_1,q_1)=1$ and $y_0$ is irrational, then
\begin{align*}
&i\left(S_{n,m,s}h(x_0,y_0);\left[0,g\left(\frac{j}{q_1}\right)\right]\right)
\geq \frac{1}{q_1},\qquad j=0,\dots,q_2-1
\end{align*}
while
\begin{align*}
&i\left(S_{n,m,s}h(x_0,y_0);\frac{1}{2} \right)=\frac{1}{q_1}\;,\quad
i\left(S_{n,m,s}h(x_0,y_0);\frac14\right)=1-\frac{1}{q_1}
\end{align*}
if $s=1$;

\item[iii)] if $x_0$ is irrational and $y_0=\frac{p_2}{q_2}$ with $p_2,q_2\in \N$, $q_2\neq 0$,
$GCD(p_2,q_2)=1$, then
\begin{align*}
&i\left(S_{n,m,s}h(x_0,y_0);\left[0,g\left(\frac{j}{q_2}\right)\right]\right)\geq\frac{1}{q_2},
\qquad j=0,\dots,q_2-1;
\end{align*}
while
\begin{align*}
&i\left(S_{n,m,s}h(x_0,y_0);\frac{1}{2} \right)=\frac{1}{q_2}\;,\quad
i\left(S_{n,m,s}h(x_0,y_0);\frac14\right)=1-\frac{1}{q_2}
\end{align*}
if $s=1$;

\item[iv)] if $x_0$ and $y_0$ are both irrational, then
\[
i\left(S_{n,m,s}h(x_0,y_0);A\right)=|G_s^{-1}(A)|
\]
for every Peano-Jordan measurable set $A\subset \R$ if $s>1$,
while
\[
i\left(S_{n,m,s}h(x_0,y_0);\frac14\right)=1
\]
if $s=1$.
\end{itemize}
\end{itemize}

%
\end{theorem}

\proof
We define the functions $h_1,\;h_2:[0,1]\to \R$ as follows
\begin{equation}\label{eq:h1-h2_inverse}
h_1(t):=\left\{
\begin{array}{ll}
1\, &\text{if }t< x_0,
\\
0\, &\text{if }t\geq x_0,
\end{array}
\right.
\qquad
h_2(t):=\left\{
\begin{array}{ll}
1\, &\text{if }t< y_0,
\\
0\, &\text{if }t\geq y_0\;,
\end{array}
\right.
\end{equation}
then we can write
\[
h(x,y)=h_1(x)h_2(y)\;,
\]
and consequently we have
\begin{align}\label{eq:decomp-Snm}
&S_{n,m}h(x,y)=S_nh_1(x)S_mh_2(y)\;,
\end{align}
moreover using \cite[Theorem 3.1]{CMT}, 
it follows that
\begin{equation}\label{eq:conv_compact1_S}
\lim_{n\rightarrow \infty}S_nh_1=h_1 \text{ uniformly on }[0,1]\setminus\{\,x_0\, \}
\end{equation}
and
\begin{equation}\label{eq:conv_compact2_S}
\lim_{m\rightarrow \infty}S_mh_2=h_2 \text{ uniformly on }[0,1]\setminus\{\,y_0\, \}\;.
\end{equation}

Let us consider a compact set $K\subset [0,1]\times[0,1]\setminus Q$; there exist
$0\leq a_1<x_0< a_2\leq 1$ and $0\leq b_1<y_0< b_2\leq1$ such that
\[
K\subset [0,a_1]\times[0,b_1]\cup [a_2,1]\times [0,1]\cup [0,a_2]\times[b_2,1].
\]
Firstly, let us consider $(x,y)\in [0,a_1]\times[0,b_1]$; then
\[
|S_{n,m}h(x,y)-h(x,y)|=|S_{n,m}h(x,y)-1|
\]
which is uniformly convergent to $0$ as $n,m\rightarrow \infty$ by \eqref{eq:conv_compact1_S}
and \eqref{eq:conv_compact2_S}. Therefore we can conclude that
\[
\lim_{n,m\rightarrow \infty}S_{n,m}h=h\qquad \text{uniformly in }[0,a_1]\times[0,b_1].
\]

Let us consider now $(x,y)\in [a_2,1]\times [0,1]$; for sufficiently large $n,m\geq 1$ there exist
$k_0,\ell_0\geq 1$ such that
\[
\frac{k_0}{n}\leq x_0<\frac{k_0+1}{n},\ \qquad \frac{\ell_0}{m}\leq y_0<\frac{\ell_0+1}{m}.
\]
Notice that
\begin{equation}\label{eq:boundness_S}
|S_{m}h_2(y)|=\frac{\sum^{\ell_0}_{j=0}|y-\frac{j}{m}|^{-s}}{\sum^m_{j=0}|y-\frac{j}{m}|^{-s}}\leq
\frac{\sum^{\ell_0}_{j=0}|y-\frac{j}{m}|^{-s}}{\sum^{\ell_0}_{j=0}|y-\frac{j}{m}|^{-s}+
\sum^m_{j=\ell_0+1}|y-\frac{j}{m}|^{-s}}\leq 1
\end{equation}
and this estimate is independent of $y\in [0,1]$ and $n\geq1$. Then
\[
|S_{n,m}h(x,y)|\leq |S_nh_1(x)|
\]
where, by \eqref{eq:conv_compact1_S}, the last term converges uniformly to 0 as $n\rightarrow \infty$. So we can conclude that
\[
\lim_{n,m\rightarrow \infty}S_{n,m}h=h\qquad \text{uniformly in } [a_2,1]\times [0,1].
\]

Arguing similarly, exchanging the role between $x$ and $y$, we get also that
\[
\lim_{n,m\rightarrow \infty}S_{n,m}h=h\qquad \text{uniformly in } [0,a_2]\times[b_2,1].
\]

The proof of the claims 1)-3) is at all similar to the one of Theorem \ref{th:lagrange-multi},
we have only to use the decomposition \eqref{eq:decomp-Snm} and \cite[Theorem 3.1]{CMT} in place of Theorem \ref{th:Lagrange}.

In particular, in the case 3)--i), if $s=1$, using \cite[Theorem 3.1]{CMT}, we can consider two subsequences $(S_{k_1(n),1}h_1(x_0))_{n\geq 1}$ and
$(S_{k_2(n),1}h_1(x_0))_{n\geq 1}$ of $(S_{n,1}h_1(x_0))_{n\geq 1}$, converging respectively to 1 and
$\frac 12$, with density respectively $\frac{1}{q_1}$ and $1-\frac{1}{q_1}$. Then
\[
\lim_{n\rightarrow \infty}S_{k_1(n),m,1}h(x_0,y_0)=S_mh_2(y_0),
\]

\[
\lim_{n\rightarrow \infty}S_{k_2(n),m,1}h(x_0,y_0)=\frac 12 S_mh_2(y_0)
\]
and the previous limits are uniform with respect to $m\geq 1$ by \eqref{eq:boundness_S}. So
we can apply Proposition \ref{pr:1-bivariate} and we get
\[
i\left(S_{n,m,1}h(x_0,y_0);\frac 12\right)\geq \frac{1}{q_1}i\left(S_{m,1}h_2(y_0);\frac 12\right)=\frac{1}{q_1}
\]
and
\[
i\left(S_{n,m,1}h(x_0,y_0);\frac 14\right)\geq \left(1-\frac{1}{q_1}\right)i\left(\frac 12 S_{m,1}h_2(y_0);
\frac 14\right)=1-\frac{1}{q_1}.
\]
Finally, the previous inequalities become equalities by Proposition \ref{pr:sum_gen_multi}.

Arguing similarly the others cases of our claim can be proved.
\qed


\begin{thebibliography}{99}



\bibitem{BDVM}
\textsc{Bojanic, R., Della Vecchia, B., and Mastroianni, G.,} \emph{On the approximation of bounded functions
with discontinuities of the first kind by generalized Shepard operators}. Acta Math. Hungarica \textbf{85} (1-2) (1999), 29--57.

\bibitem{CMT}
\textsc{Campiti, M., Mazzone, G., and Tacelli, C.,} \emph{Interpolation of discontinuous functions}. Preprint, 2010, available at http://arxiv.org, arXiv:1009.2954v2.

\bibitem{C1}
\textsc{Connor, J. S.,} \emph{The statistical and strong $p$-Ces\`aro convergence of sequences}. Analysis \textbf{8} (1988), 47--63.

\bibitem{DVM}
\textsc{Della Vecchia, B. and Mastroianni, G.,} \emph{On functions approximation by Shepard-type operators -- a survey},
in Approximation Theory, Wavelets and Applications (Maratea, 1994), 335--346.
NATO Adv. Sci. Inst. Ser. C Math. Phys. Sci. \textbf{454}, Kluwer Acad. Publ., Dordrecht, 1995.

\bibitem{DO}
\textsc{Duman, O. and Orhan, C.,} \emph{Statistical approximation by positive linear operators}.
Studia Math. \textbf{161} (2) (2004), 187--197.

\bibitem{F1}
\textsc{Fridy, J. A.,} \emph{On statistical convergence}. Analysis \textbf{5} (1985), 301--313.

\bibitem{F2}
\textsc{Fridy, J. A.,} \emph{Statistical limit points}. Proc. Amer. Math. Soc. \textbf{118} (1993), 1187--1192.

\bibitem{F3}
\textsc{Fridy, J. A. and Orhan, C.,} \emph{Statistical limit superior and limit inferior}. Proc. Amer. Math.
Soc. \textbf{125} (1997), 3625--3631.

\bibitem{G}
\textsc{Gal, S. G.,} Global smoothness and shape preserving interpolation by classical operators.
\emph{Birk\"auser} (2005).

\bibitem{MM}
\textsc{Mastroianni, G. and Milovanovi\'c, G.,} Interpolation processes - Basic theory and applications.
\emph{Springer Monographs in Mathematics, Springer-Verlag Berlin Heidelberg}, 2008.

\bibitem{Sh}
\textsc{Shepard, D.,} \emph{A two-dimensional interpolation function for irregularly spaced points}.
Proceedings 1968 Assoc. Comput. Machinery National Conference, 517--524.

\bibitem{S}
\textsc{Szeg\"o, G., Orthogonal Polynomials.} \emph{Colloquium Publ., vol. XXIII, Amer. Math. Soc.},
Providence, RI, 1959; Russian translation: Fizmatlit, Moscow, 1962.

\bibitem{V}
\textsc{Vertesi, P.,} \emph{Lagrange interpolation for continuous functions of bounded variation}.
Acta Math. Acad. Scient. Hung. \textbf{35} (1-2) (1980), 23--31.

\bibitem{WN}
\textsc{Williams, K. S. and Nan-Yue, Z.,} \emph{Special values of the Lerch zeta function and the evaluation of certain integrals}.
Proc. Amer. Math. Soc. \textbf{119} (1) (1993), 35--49.

\bibitem{WW}
\textsc{Whittaker, E. T. and Watson, G. N.,} A course of modern analysis, 4th ed. \emph{Cambridge Univ. Press},
Cambridge and New York, (1963).


\end{thebibliography}
\end{document}